\documentclass{amsart}

\usepackage{amssymb,amsmath,latexsym,amsfonts,amsthm,mathrsfs}
\usepackage{enumerate}
\newtheorem{thm}{Theorem}
\newtheorem{defn}{Definition}

\newtheorem*{claim}{Claim}

\newtheorem{lemma}{Lemma}

\newtheorem{cor}{Corollary}

\newcommand{\R}{\mathbb{R}}
\newcommand{\Q}{\mathbb{Q}}
\newcommand{\funct}[2]{#1 \longrightarrow #2}

\newcommand{\Ur}{\textbf{Q}_{S}}
\newcommand{\arrows}[3]{\longrightarrow {#1}^{#2}_{#3}}
\newcommand{\B}{\mathcal{B} _S}
\newcommand{\U}{\mathcal{U}_S}
\newcommand{\m}[1]{\textbf{#1}}
\newcommand{\mc}[1]{\widetilde{\textbf{#1}}}
\newcommand{\restrict}[2]{#1 \upharpoonright #2}

\author{L. Nguyen Van Th\'e}

\address{Equipe de Logique Math\'ematique, UFR de Math\'ematiques,
 (case 7012), Universit\'e Denis Diderot Paris 7, 2 Place Jussieu,
 75251 Paris Cedex 05, France.}

\email{nguyenl@logique.jussieu.fr}

\title{Big Ramsey Degrees and Divisibility in Classes of Ultrametric Spaces}

\date{July, 2005}

\begin{document}

\begin{abstract}

Given a countable set $S$ of positive reals, we study finite-dimensional Ramsey-theoretic properties of the countable ultrametric Urysohn space $\textbf{Q} _S$ with distances in $S$.

\end{abstract}

\maketitle

\section{Introduction}

This note is the continuation of our paper \cite{NVT} where Ramsey-type properties of several classes of finite ultrametric spaces are studied. Recall that a metric space $\m{X} = (X, d^\m{X})$ is \emph{ultrametric} when given any $x, y, z$ in $\m{X}$,
\[d^\m{X}(x,z) \leqslant \max(d^\m{X}(x,y), d^\m{X}(y,z))\]

Given $S \subset ]0, + \infty[$, the class $\U$ is defined as the class of all finite ultrametric spaces with strictly positive distances in $S$. It turns out that when $S$ is at most countable, there is, up to isometry, a unique countable metric space $\Ur$ such that i) The family of finite metric subspaces of $\Ur$ is exactly $\U$ and ii) $\Ur$ is \emph{ultrahomogeneous}, that is any isometry between any two finite subspaces of $\Ur$ can be extended to an isometry of $\Ur$ onto itself. $\Ur$ is called \emph{Urysohn space associated to $\U$}. It is a variation of the rational Urysohn space $\textbf{U} _0$ constructed by Urysohn in \cite{U}: The difference between $\Ur$ and $\textbf{U} _0$ is that whereas $\Ur$ is the countable universal ultrahomogeneous space attached to $\U$, $\textbf{U} _0$ is related to the class $\mathcal{M} _{\Q}$ of all finite metric spaces with rational distances. Unlike $\textbf{U} _0$, $\Ur$ can be represented quite simply. Namely, $\Ur$ can be seen
  as the set of all finitely supported elements of $\omega ^S$ equipped with the distance $d^{\Ur}$ defined by 
\[ d^{\Ur}(x,y) = \max \{s \in S : x(s) \neq y(s) \} \]

The role that $\textbf{U} _0$ and $\Ur$ play with respect to $\mathcal{M} _{\Q}$ and $\U$ respectively are exactly the same as the role that the Rado graph $\mathcal{R}$ plays for the class of finite graphs. Our concern here is to obtain for $\Ur$ the analogs of well-known results of the form 
\[ \mathcal{R} \arrows{\mathcal{(R)}}{\m{G}}{k, l} \]

More precisely, for metric spaces spaces $\m{X}$, $\m{Y}$ and $\m{Z}$, write $\m{X} \cong \m{Y}$ when there is an isometry from $\m{X}$ onto $\m{Y}$ and define the set $\binom{\m{Z}}{\m{X}}$ as 
\[ \binom{\m{Z}}{\m{X}} = \{ \mc{X} \subset \m{Z} : \mc{X} \cong \m{X} \} \]

Then, given positive integers $k,l$, the symbol $ \m{Z} \arrows{(\m{Y})}{\m{X}}{k,l} $ abbreviates the fact that:
\begin{center} For any $\chi : \funct{\binom{\m{Z}}{\m{X}}}{k}$ there is $\widetilde{\m{Y}} \in \binom{\m{Z}}{\m{Y}}$ such that $\chi$ does not take more than $l$ values on $\binom{\widetilde{\m{Y}}}{\m{X}}$. \end{center} 

In \cite{NVT}, it was proved that any element $\m{X}$ has a \emph{finite Ramsey degree in $\U$}, meaning that there is an integer $l$ (depending on $\m{X}$) for which given any $\m{Y} \in \U$ and any $k \in \omega \smallsetminus \{ 0 \}$,
there exists $\m{Z} \in \mathcal{K} $ such that 
\[ \m{Z} \arrows{(\m{Y})}{\m{X}}{k,l} \]

Consequently, 
\[ \forall k \in \omega \smallsetminus \{ 0 \} \ \ \Ur \arrows{(\m{Y})}{\m{X}}{k,l} \]

The purpose of this paper is to present the conditions on $S \subset ]0, + \infty[$ and $\m{X} \in \U$ under which this latter result remains valid when $\m{Y}$ is replaced by $\Ur$ and if so to compute the least integer $l$ such that 
\[ \forall k \in \omega \smallsetminus \{ 0 \} \ \ \Ur \arrows{(\Ur)}{\m{X}}{k,l} \]

When defined, this least $l$ is called the \emph{big Ramsey degree of $\m{X}$ in $\U$} and is in fact part of the more general notion of big Ramsey degree for an arbitrary class of finite structures, a concept defined in \cite{KPT} by Kechris, Pestov and Todorcevic in the general setting of oscillation stability of topological groups.

\begin{thm}

\label{thm:3}

Let $S$ be a finite subset of $]0, + \infty[$. Then every element of $\U$ has a big Ramsey degree in $\U$.

\end{thm}

\begin{thm}

\label{thm:4}

Let $S$ be an infinite countable subset of $]0, + \infty[$ and let $\m{X}$ be in $\U$ such that $|\m{X}| \geqslant 2$. Then $\m{X}$ does not have a big Ramsey degree in $\U$.
  
\end{thm}

Theorem \ref{thm:4} does not cover the case $|X|=1$, which is related to the \emph{divisibility} properties of $\Ur$. A metric space $\m{X}$ is \emph{indivisible} when given any $k \in \omega \smallsetminus \{ 0 \}$ and any map $\chi : \funct{\m{X}}{k}$, there is an isometric copy $\mc{X}$ of $\m{X}$ included in $\m{X}$ on which $\chi$ is constant. Otherwise, $\m{X}$ is \emph{divisible}. Our results read as follows:

\begin{thm}

\label{thm:1}

Let $S$ be an infinite countable subset of $]0, + \infty[$ and assume that the reverse linear ordering $>$ on $\R$ does not induce a well-ordering on $S$. Then there is a map $\chi : \funct{\Ur}{\omega}$ whose restriction on any isometric copy $X$ of $\Ur$ inside $\Ur$ has range $\omega$. 
\end{thm} 

In particular, in this case, $\Ur$ is divisible.

\begin{thm}

\label{thm:2}

Let $S$ be an infinite countable subset of $]0, + \infty[$ and assume that the reverse linear ordering $>$ on $\R$ induces a well-ordering on $S$. Then $\Ur$ is indivisible.
\end{thm}

It should be mentionned at that point that we are now aware of the fact that theorem \ref{thm:1} and theorem \ref{thm:2} were first obtained completely independently by Delhomm\'e, Laflamme, Pouzet and Sauer in \cite{DLPS} where a precise analysis of divisibility is carried out in the realm of countable metric spaces. In particular, \cite{DLPS} provides a necessary condition for indivisibility and solves the indivisibility problem for several countable Urysohn spaces. For example, every sphere of $\textbf{U} _0$ is divisible. 

We finish with a consequence of theorem \ref{thm:2}:

\begin{thm}

\label{thm:5}

Let $S$ be an infinite countable subset of $]0, + \infty[$ and assume that the reverse linear ordering $>$ on $\R$ induces a well-ordering on $S$. Then given any map $f : \funct{\Ur}{\omega}$, there is an isometric copy $X$ of $\Ur$ inside $\Ur$ such that $f$ is continuous or injective on $X$. 

\end{thm}

The paper is organized as follows: In the forthcoming section, we recall the connection between the notions of trees and ultrametric spaces, and use classical Ramsey's theorem to prove theorem \ref{thm:3} and theorem \ref{thm:4}. In section 3, we concentrate on the divisibility properties of $\Ur$, and prove theorem \ref{thm:1} and theorem \ref{thm:2}. Finally, we close in section 4 with the proof of theorem \ref{thm:5}. 

I would like to thank sincerely Stevo Todorcevic for his various suggestions and comments concerning this paper. I am also indebted to Jordi Lopez Abad for the numerous ideas his constant help, support and enthusiasm brought.

\section{Big Ramsey degrees in $\U$}

The purpose of this section is to provide the proofs for theorem \ref{thm:3} and theorem \ref{thm:4}. The ideas we use to reach this goal are not new. The way we met them is through some unpublished work of Galvin, but in \cite{M}, Milner writes that they were also known to and exploited by several other authors, among whom Hajnal (who apparently realized first the equivalent of lemma \ref{lemma:3} and stated it explicitly in \cite{H}), and Haddad and Sabbagh (\cite{HS1}, \cite{HS2} and \cite{HS3}).

Since everything here is connected to the notion of tree, let us start with some general facts about these objects. 

A \emph{tree} $\m{T} = (T,<^{\m{T}}) $ is a partially ordered set such that given any element $t
\in T$, the set $\{ s \in T : s <^{\m{T}} t \}$ is
$<^{\m{T}}$-well-ordered. $e(\m{T})$ then denotes the set of all linear orderings on $T$ which extend $<^{\m{T}}$.
When $\m{U}$ is also a tree, $\m{T}$ and $\m{U}$ are \emph{isomorphic} when there is an order-preserving bijection from $\m{T}$ to $\m{U}$ (in symbols, $\m{T} \cong \m{U}$).
When every element of $T$ has finitely many $<^{\m{T}}$-predecessors, the \emph{height of} $t
\in \m{T}$ is \[ \mathrm{ht}(t) = |\{ s \in T : s <^{\m{T}} t \}| \]

When $S$ is finite and given by elements $s_0 > s_1 \ldots > s_{|S|-1} > 0$, it will be convenient to see the space $\Ur$ as the set $\omega^{|S|}$ of maximal nodes of the tree $\omega^{\leqslant |S|} = \bigcup_{i\leqslant n} \omega^i$ ordered by set-theoretic inclusion and equipped with the metric defined for $x \neq y$ by \[d(x,y) = s_{\Delta (x,y)}\]

where $\Delta (x,y)$ is the height of the largest common predecessor of $x$ and $y$ in $\omega^{\leqslant |S|}$. For $A \subset \omega ^{|S|}$, set \[ A^\downarrow = \{ \restrict{a}{k} : a \in A \wedge k \leqslant n \} \]

It should be clear that when $A, B \subset \omega ^{|S|}$, then $A$ and $B$ are isometric iff $A^\downarrow \cong B^\downarrow$. Consequently, when $\m{X} \in \U$, one can define the natural tree associated to $\m{X}$ in $\U$ to be the unique (up to isomorphism) subtree $\m{T}_{\m{X}}$ of $\omega^{\leqslant |S|}$ such that for any copy $\mc{X}$ of $\m{X}$ in $\Ur$, $\mc{X} ^\downarrow \cong \m{T}_{\m{X}}$.

We now introduce some notations for the partition calculus on trees.
Given a subtree $\m{T}$ of $\omega ^{|S|}$, set \[ \binom{\omega ^{\leqslant|S|}}{\m{T}} = \{ \mc{T} : \mc{T} \subset \omega ^{\leqslant|S|} \wedge \mc{T} \cong \m{T} \} \]

When $k, l \in \omega \smallsetminus \{ 0 \}$ and for any $\chi : \funct{\binom{\omega ^{\leqslant|S|}}{\m{T}}}{k}$ there is
$\m{U} \in \binom{\omega ^{\leqslant|S|}}{\omega ^{\leqslant|S|}}$ such that $\chi$ takes at most $l$ values on $\binom{\m{U}}{\m{T}}$, we write \[ \omega ^{\leqslant|S|}
\arrows{(\omega ^{\leqslant|S|})}{\m{T}}{k,l} \]

If there is $l \in \omega \smallsetminus \{ 0 \}$ such that for any $k \in \omega \smallsetminus \{ 0 \}$, $\omega ^{\leqslant|S|}
\arrows{(\omega ^{\leqslant|S|})}{\m{T}}{k,l}$, the least such $l$ is called the \emph{Ramsey degree} of $\m{T}$ in $\omega ^{\leqslant|S|}$. 

\begin{lemma}

\label{lemma:3}

Let $X \subset \omega ^{|S|}$ and let $\m{T} = X^\downarrow$. Then $\m{T}$ has a Ramsey degree in $\omega ^{\leqslant|S|}$ equal to $|e(\m{T})|$. 
\end{lemma}

\begin{proof}
Say that a subtree $\m{U}$ of $\omega ^{\leqslant|S|}$ is \emph{expanded} when: 

i) Elements of $\m{U}$ are strictly increasing. 

ii) For every $u, v \in \m{U}$ and every $k \in |S|$, \[u(k) \neq v(k) \rightarrow (\forall j \geqslant k \ \ u(j) \neq v(j)) \]

Note that every expanded $\mc{T} \in \binom{\omega ^{\leqslant|S|}}{\m{T}}$ is linearly ordered by $\prec^{\mc{T}}$ defined by \begin{center} $s \prec^{\mc{T}} t$ iff ($s=\emptyset$ or $s(|s|) < t(|t|)$) \end{center} 

and that then $\prec^{\mc{T}}$ is a linear extension of the tree ordering on $\mc{T}$.

Now, given $\prec \in e(\m{T})$, let $\binom{\omega ^{\leqslant|S|}}{\m{T},\prec}$ denote the set of all expanded $\mc{T} \in \binom{\omega ^{\leqslant|S|}}{\m{T}}$ \emph{with type} $\prec$, that is such that the order-preserving bijection between the linear orderings $(\mc{T}, \prec^{\mc{T}})$ and $(\m{T}, \prec)$ induces an isomorphism between the trees $\mc{T}$ and $\m{T}$. Define the map $\psi _{\prec} : \funct{\binom{\omega ^{\leqslant|S|}}{\m{T},\prec}}{[\omega]^{|\m{T}|-1}}$ by 
\[ \psi _{\prec} (\mc{T}) = \{ t(|t|) : t \in \mc{T} \smallsetminus \{ \emptyset \} \} \] 

Then $\psi _{\prec}$ is a bijection. Call $\varphi _{\prec}$ its inverse map. 

Now, let $k \in \omega \smallsetminus \{ 0 \}$ and $\chi : \funct{\binom{\omega ^{\leqslant|S|}}{\m{T}}}{k}$. Define $\Lambda : \funct{[\omega]^{|T|-1}}{k^{e(\m{T})}}$ by 
\[ \Lambda (M) = (\chi (\varphi _{\prec}(M)))_{\prec \in e(\m{T})} \]  

By Ramsey's theorem, find an infinite $N \subset \omega$ such that $\Lambda$ is constant on $[N]^{|\m{T}|-1}$. Then, on the subtree $N^{\leqslant |S|}$ of $\omega^{\leqslant |S|}$, any two expanded elements of $\binom{\omega ^{\leqslant|S|}}{\m{T}}$ with same type have the same $\chi$-color. Now, let $\m{U}$ be an expanded everywhere infinitely branching subtree of $N^{\leqslant |S|}$. Then $\m{U}$ is isomorphic to $\omega^{\leqslant |S|}$ and $\chi$ does not take more than $|e(\m{T})|$ values on $\binom{\m{U}}{\m{T}}$.

To finish the proof, it remains to show that $|e(\m{T})|$ is the best possible bound. To do that, simply observe that for any $\m{U} \in \binom{\omega^{\leqslant |S|}}{\omega^{\leqslant |S|}}$, every possible type appears on $\binom{\m{U}}{\m{T}}$.

\end{proof}

This lemma has two direct consequences concerning the existence of big Ramsey degrees in $\U$. Indeed, it should be clear that when $\m{X} \in \U$, $\m{X}$ has a big Ramsey degree in $\U$ iff $\m{T}_{\m{X}}$ has a Ramsey degree in $\omega ^{\leqslant|S|}$ and that these degrees are equal. So on the one hand:

\begin{cor}[Theorem \ref{thm:3}]

Let $S$ be a finite subset of $]0, + \infty[$. Then every element of $\U$ has a big Ramsey degree in $\U$.

\end{cor} 

On the other hand, observe that if $S \subsetneq S'$ are finite and $\m{X} \in \U$ has size at least two, then the big Ramsey degree $T_{\mathcal{U} _{S'}}(\m{X})$ of $\m{X}$ in $\mathcal{U} _{S'}$ is strictly larger than the big Ramsey degree of $\m{X}$ in $\U$. In particular, $T_{\mathcal{U} _{S'}}(\m{X})$ tends to infinity when $|S'|$ tends to infinity. This fact has the following consequence:

\begin{cor}[Theorem \ref{thm:4}]

Let $S$ be an infinite countable subset of $]0, + \infty[$ and let $\m{X}$ be in $\U$ such that $|\m{X}| \geqslant 2$. Then $\m{X}$ does not have a big Ramsey degree in $\U$.

\end{cor} 

\begin{proof}

It suffices to show that for every $k \in \omega \smallsetminus \{ 0 \}$, there is $k' > k$ and a coloring $\chi : \funct{\binom{\Ur}{\m{X}}}{k'}$ such that for every $Q \in \binom{\Ur}{\Ur}$, the restriction of $\chi$ on $\binom{Q}{\m{X}}$ has range $k'$. 

Thanks to the previous remark, we can fix $S' \subset S$ finite such that $X \in \mathcal{U}_{S'}$ and the big Ramsey degree $k'$ of $\m{X}$ in $\mathcal{U}_{S'}$ is larger than $k$. Recall that $\Ur \subset \omega ^S$ so if $\textbf{1} _{S'} : \funct{S}{2}$ is the characteristic function of $S'$, it makes sense to define $f : \funct{\Ur}{\textbf{Q}_{S'}}$ by 
\[ f(x) = \textbf{1}_{S'} x\]

Observe that $d(f(x),f(y)) = d(x,y)$ whenever $d(x,y) \in S'$. Thus, given any $Q \in \binom{\Ur}{\Ur}$, the direct image $f''Q$ of $Q$ under $f$ is in $\binom{\textbf{Q}_{S'}}{\textbf{Q}_{S'}}$. Now, let $\chi ' : \funct{\binom{\textbf{Q}_{S'}}{\m{X}}}{k'}$ be such that for every $Q' \in \binom{\textbf{Q}_{S'}}{\textbf{Q}_{S'}}$, the restriction of $\chi '$ to $\binom{Q'}{\m{X}}$ has range $k'$. Then $\chi = \chi ' \circ f $ is as required. 

\end{proof}

\section{Divisibility properties of $\Ur$}

In this section, we study the divisibility properties of $\Ur$ and provide the proofs for theorem \ref{thm:1} and theorem \ref{thm:2}. Recall that a metric space $\m{X}$ is \emph{indivisible} when given any $ k \in \omega$ and any map $\chi : \funct{\m{X}}{k}$, there is an isometric copy $\mc{X}$ of $\m{X}$ on which $\chi$ is constant. Otherwise, $\m{X}$ is \emph{divisible}. 

Unlike theorem \ref{thm:3} and theorem \ref{thm:4}, the proofs here do not use any classical partition result via a particular coding of the objects, but rather take place on the geometrical level.    

For notational convenience, we will often simply write $d$ instead of $d^{\Ur}$.

\subsection{Proof of theorem \ref{thm:1}}

Fix an infinite and countable subset $S$ of $]0, + \infty[$ such that the reverse linear ordering $>$ on $\R$ does not induce a well-ordering on $S$. The idea to prove that $\Ur$ is divisible is to use a coloring which is constant on some particular spheres. 

More precisely, observe that $(S, >)$ not being well-ordered, there is a strictly increasing sequence $(s_i)_{i \in \omega}$ of reals such that $s_0 = 0$ and $s_i \in S$ for every $i > 0$. Observe that we can construct a subset $E$ of $\Ur$ such that given any $y \in \Ur$, there is exactly one $x$ in $E$ such that for some $i < \omega$, $d(x,y) < s_i$. Indeed, if $\sup _{i < \omega}s_i = \infty$, simply take $E$ to be any singleton. Otherwise, let $\rho = \sup _{i < \omega}s_i$ and choose $E \subset \Ur$ maximal such that 
\[ \forall x, y \in E \ \ d(x,y) \geqslant \rho \]

To define $\chi : \funct{\Ur}{\omega}$, let $(A_j)_{j \in \omega}$ be a family of infinite pairwise disjoint subsets of $\omega$ whose union is $\omega$. Then, for $y \in \Ur$, let $e(y)$ and $i(y)$ be the unique elements of $E$ and $\omega$ respectively such that $d(e(y),y) \in [s_{i(y)} , s_{i(y)+1}[$, and set \begin{center} $\chi (y) = j$ iff $i(y) \in A_j$ \end{center} 

\begin{claim}
$\chi$ is as required.
\end{claim}

\begin{proof}
Let $Y \subset \Ur$ be isometric to $\Ur$. Fix $y \in Y$. For every $j \in \omega$, pick $i_j > i(y) + 1$ such that $i_j \in A_j$. Since $Y$ is isometric to $\Ur$, we can find an element $y_j$ in $Y$ such that $d(y, y_j) = s_{i_j}$. We claim that $\chi (y_j) = j$, or equivalently $i(y_j) \in A_j$. Indeed, consider the triangle $\{e(y), y, y_j \}$. Observe that in an ultrametric space every triangle is isosceles with short base and that here, 
\[ d(e(y),y) < s_{i_j} = d(y,y_j) \]

Thus, 
\[ d(e(y),y_j) = d(y,y_j) \in [s_{i_j},s_{i_j +1}[ \]

And therefore $e(y_j) = e(y)$ and $i(y_j) = i_j \in A_j$. 

\end{proof}

\subsection{Proof of theorem \ref{thm:2}}
Fix an infinite countable subset $S$ of $]0, + \infty[$ such that the reverse linear ordering $>$ on $\R$ induces a well-ordering on $S$. Our goal here is to show that the space $\Ur$ is indivisible. 

Observe first that the collection $\B$ of metric balls of $\Ur$ is a
tree when ordered by reverse set-theoretic inclusion. When $x \in \Ur$
and $r \in S$, $B(x,r)$ denotes the set $\{ y \in \Ur : d^{\Ur}(x,y)
\leqslant r \}$. $x$ is called a \emph{center} of the ball and $r$ a
\emph{radius}. Note that in $\Ur$, non empty balls have a unique radius but
admit all of their elements as centers. Note also that when $s > 0 $
is in $S$, the fact that $(S,>)$ is well ordered  allows to define 
\[ s^- = \max \{t \in S : t < s \} \]

The main ingredients are contained in the following definition and lemma. 

\begin{defn}

\label{defn:1}

Let $A \subset \Ur$ and $b \in \B$ with radius $r \in S \cup \{ 0 \}$. Say that \emph{$A$ is
small in $b$} when $r=0$ and $A \cap b = \emptyset$ or $r > 0$
and $A \cap b$ can be covered by finitely many balls of radius $r^-$.
\end{defn}

We start with an observation. Assume that $\{ x_n : n \in \omega \}$ is an enumeration of $\Ur$, and that we are trying to build inductively a copy $\{ a_n : n \in \omega \}$ of $\Ur$ in $A$ such that for every $n, m \in \omega$, $d(a_n , a_m) = d(x_n , x_m)$. Then the fact that we may be blocked at some finite stage exactly means that at that stage, a particular metric ball $b$ with $A \cap b \neq \emptyset$ is such that $A$ is small in $b$. This idea is expressed in the following lemma.

\begin{lemma}

\label{lemma:1}

Let $X \subset \Ur$. TFAE :

i) $\binom{X}{\Ur} \neq \emptyset$.

ii) There is $Y \subset X$ such that $Y$ is not small in $b$ whenever $b \in \B$ and $Y \cap b \neq \emptyset$.

\end{lemma}

\begin{proof}
Assume that i) holds and let $Y$ be a copy of $\Ur$ in $X$. Fix $b \in \B$ with radius $r$ and such that $Y \cap b \neq \emptyset$. Pick $x \in Y \cap b $ and let $E \subset \Ur$ be an infinite subset where all the distances are equal to $r$. Since $Y$ is isometric to $\Ur$, $Y$ includes a copy $\tilde{E}$ of $E$ such that $x \in \tilde{E}$. Then $\tilde{E} \subset Y \cap b$ and cannot be covered by finitely many balls of radius $r^-$, so ii) holds.

Conversely, assume that ii) holds. Let $\{x_n : n \in \omega \}$ be an enumeration of the elements of $\Ur$. We are going to construct inductively a sequence $(y_n)_{n \in \omega}$ of elements of $Y$ such that \[ \forall m, n \in \omega \ \ d(y_m , y_n) = d(x_m , x_n) \]

For $y_0$, take any element in $Y$. In general, if $(y_n)_{n \leqslant k}$ is built, construct $y_{k+1}$ as follows. Consider the set $E$ defined as 
\[ E = \{ y \in \Ur : \forall \ n \leqslant k \ \ d(y,y_n) = d(x _{k+1} , x_n) \} \]

Let also 
\[ r = \min \{ d(x_{k+1},x_n) : n \leqslant k \} \]

and
\[ M = \{ n \leqslant k : d(x_{k+1},x_n) = r \} \] 

We want to show that $E \cap Y \neq \emptyset$. Observe first that for every $m, n \in M$, $d(y_m , y_n) \leqslant r$. Indeed,
\[d(y_m , y_n ) =  d(x_m , x_n) \leqslant \max (d(x_m , x_{k+1}) , d(x_{k+1} , x_n)) = r \]

So in particular, all the elements of $\{ y_m : m \in M \}$ are contained in the same ball $b$ of radius $r$.

\begin{claim}
$E = b \smallsetminus \bigcup _{m \in M} B(y_m, r^-)$.
\end{claim}

\begin{proof}
It should be clear that \[ E \subset b \smallsetminus \bigcup _{m \in M} B(y_m, r^-) \] 

On the other hand, let $y \in b \smallsetminus \bigcup _{m \in M} B(y_m, r^-)$. Then for every $m \in M$, \[ d(y,y_m) = r = d(x_{k+1}, x_m) \]

so it remains to show that $d(y,y_n) = d(x_{k+1}, x_n)$ whenever $n \notin M$. To do that, we use again the fact that every triangle is isosceles with short base. Let $m \in M$. In the triangle $\{x_m, x_n, x_{k+1} \}$, we have $d(x_{k+1}, x_n) > r$ so \[ d(x_m, x_{k+1}) = r < d(x_n , x_m) = d(x_n , x_{k+1}) \]

Now, in the triangle $\{ y_m, y_n, y\}$, $d(y, y_m) = r$ and  $d(y_m , y_n) = d(x_m, x_n) > r$. Therefore, \[ d(y, y_n) = d(y_m, y_n) = d(x_m, x_n) = d(x_{k+1}, x_n) \ \  \] 

\end{proof}

We consequently need to show that $(b \smallsetminus \bigcup _{m \in M} B(y_m, r^-)) \cap Y \neq \emptyset$. To achieve that, simply observe that when $m \in M$, we have $y_m \in Y \cap b$. Thus, $Y \cap b \neq \emptyset$ and by property ii), $Y$ is not small in $b$. In particular, $Y \cap b $ is not included in $\bigcup _{m \in M} B(y_m, r^-)$. 

\end{proof}

We are now ready to prove theorem \ref{thm:2}. However, before we do so, let us make another observation concerning the notion smallness. Let $\Ur = A \cup B$. 

Note that if $A$ is small in $b \in \B$, then 1) $A \cap b$ cannot contribute to build a copy of $\Ur$ in $A$ and 2) $B \cap b$ is isometric to $b$. So intuitively, everything happens as if $b$ were completely included in $B$. So the idea is to remove from $A$ all those parts which are not essential and to see what is left at the end. More precisely, define a sequence $(A _{\alpha})_{\alpha \in \omega _1}$ recursively as follows: 

\begin{itemize}

\item $A_0 = A$.

\item $A_{\alpha + 1} = A_{\alpha} \smallsetminus \bigcup \{ b : A_{\alpha} \ \mathrm{is \ small \ in \ b} \}$. 

\item For $\alpha < \omega _1$ limit, $A_{\alpha} = \bigcap _{\eta < \alpha} A_{\eta}$. 

\end{itemize}

Since $\Ur$ is countable, the sequence is eventually constant. Set \begin{center} $\beta = \min \{ \alpha < \omega _1 : A_{\alpha +1} = A_{\alpha} \}$ \end{center}

Observe that if $A_{\beta}$ is non-empty, then $A_{\beta}$ is not small in any metric ball it intersects. Indeed, suppose that $b \in \B$ is such that $A_{\beta}$ is small in $b$. Then $A_{\beta + 1} \cap b = \emptyset$. But $A_{\beta + 1} = A_{\beta}$ so $A_{\beta} \cap b = \emptyset$. Therefore, since $A_{\beta} \subset A$, $A$ satisfies condition ii) of lemma \ref{lemma:1} and $\binom{A}{\Ur} \neq \emptyset$.

It remains to consider the case where $A_{\beta} =
\emptyset$. According to our second observation, the intuition is that
$A$ is then unable to carry any copy of $\Ur$ and is only composed of
parts which do not affect the metric structure of $B$. Thus, $B$
should include an isometric copy of $\Ur$. For $\alpha < \omega _1$,
let $\mathcal{C}_{\alpha}$ be the set of all minimal elements (in the
sense of the tree structure on $\B$) of the collection $\{ b \in \B :
A_{\alpha} \ \mathrm{is \ small \ in \ b} \}$. Note that since all
points of $B$ can be seen as balls of radius $0$ in which $A$ is
small, we have $B \subset \bigcup \mathcal{C}_0$. Note also that
$(\bigcup \mathcal{C}_{\alpha})_{\alpha < \omega _1}$ is increasing. By induction on $\alpha > 0$, it
follows that
\[ \forall \ 0 < \alpha < \omega _1 \ \ A_{\alpha} = \Ur \smallsetminus
\bigcup_{\eta < \alpha} \bigcup \mathcal{C}_{\eta} \ \ \ \ (*)\] 

\begin{claim}
Let $\alpha < \omega _1$ and $b \in \mathcal{C}_{\alpha}$ with radius $r \in S$. Then there are $c_0 \ldots c_{n-1}$ in $\B$ with radius $r^-$ and included in $b$ such that \[ b = \bigcup _{i<n} c_i \cup \bigcup _{\eta < \alpha} \bigcup \{ c \in \mathcal{C}_{\eta} : c \subset b \} \]

\end{claim}

\begin{proof}
$A_{\alpha}$ is small in b so find $c_0 \ldots c_{n-1} \in \B$ with radius $r^-$ and included in $b$ such that \[ A_{\alpha} \cap b \subset \bigcup_{i<n}c_i \]

Then thanks to $(*)$ \[ b \smallsetminus \bigcup_{i<n}c_i \subset \bigcup_{\eta < \alpha} \bigcup \mathcal{C}_{\eta} \]

Note that by minimality of $b$, if $\eta < \alpha$, then $b \subsetneq c $ cannot happen for any element of $\mathcal{C}_{\eta}$. It follows that either $c \cap b = \emptyset$ or $c \subset b$. Therefore, \[ b \smallsetminus \bigcup_{i<n}c_i \subset \bigcup_{\eta < \alpha} \bigcup \{ c \in \mathcal{C}_{\eta} : c \subset b\} \ \  \]  

\end{proof}

\begin{claim}
Let $\alpha < \omega _1$ and $b \in \mathcal{C}_{\alpha}$. Then $\binom{B \cap b}{b} \neq \emptyset$. 
\end{claim}

\begin{proof}
We proceed by induction on $\alpha < \omega _1$. 

For $\alpha = 0$, let $b \in \mathcal{C}_0$. Without loss of generality, we may assume that the radius $r$ of $b$ is strictly positive and hence in $S$. $A_0 = A$ is small in $b$ so find $c_0, \ldots , c_{n-1}$ with
radius $r^-$ such that $A \cap b \subset \bigcup_{i<n}c_i$. Then $b
\smallsetminus \bigcup_{i<n}c_i$ is isometric to $b$ and is included
in $B \cap b$.

Suppose now that the claim is true for every $\eta < \alpha$. Let $b \in \mathcal{C}_{\alpha}$ with radius $r \in S$. Thanks to the previous claim, we can find $c_0 \ldots c_{n-1} \in \B$ with radius $r^-$ and included in $b$ such that \[
b = \bigcup _{i<n} c_i \cup \bigcup _{\eta < \alpha} \bigcup \{ c \in \mathcal{C}_{\eta} : c \subset b \}
\]

Observe that \[ \bigcup _{\eta < \alpha} \bigcup \{ c \in \mathcal{C}_{\eta} : c \subset b \} = \bigcup \{ c \in \bigcup _{\eta < \alpha} : c \subset b \} \]

It follows that if $\mathcal{D}_{\alpha}$ is defined as the set of all minimal
elements (still in the sense of the tree structure on $\B$) of the
collection \[ \{ c \in \bigcup _{\eta < \alpha} \mathcal{C}_{\eta} : c \subset b \wedge \forall i < n \ \ c \cap c_i = \emptyset \} \]

Then $\{c_i : i < n \} \cup \mathcal{D}_{\alpha}$ is a collection of pairwise disjoint balls and $\bigcup \mathcal{D}_{\alpha}$ is isometric to $b$. By induction hypothesis, $\binom{B \cap c}{c} \neq \emptyset$ whenever $c \in \mathcal{D}_{\alpha}$ and there is an isometry $\varphi _c : \funct{c}{B \cap c}$. Now, let $\varphi : \funct{\bigcup \mathcal{D}_{\alpha}}{B \cap b}$ be defined as \[ \varphi = \bigcup_{c \in \mathcal{D}_{\alpha}} \varphi _c \] 

We claim that $\varphi$ is an isometry. Indeed, let $x, x' \in \bigcup \mathcal{D}_{\alpha}$. If there is $c \in \mathcal{D}_{\alpha}$ such that $x, x' \in c$ then \[ d(\varphi (x) , \varphi (x')) = d(\varphi _c (x) , \varphi _c (x')) = d(x,x') \]

Otherwise, find $c \neq c' \in \mathcal{D}_{\alpha}$ with $x \in c$ and $x' \in c'$. Observe that since we are in an ultrametric space, we have \[ \forall y, z \in c \ \ \forall y', z' \in c' \ \ d(y,y') = d(z,z') \] 

Thus, since $x, \varphi (x) \in c$ and $x', \varphi(x') \in c'$, we get \[ d(\varphi (x) , \varphi (x')) = d(x,x') \] 

\end{proof}

To finish the proof of the theorem, it suffices to notice that as a metric ball (the unique ball of radius $\max S$), $\Ur$ is in $\mathcal{C}_{\beta}$. So according to the previous claim, $\binom{B}{\Ur} \neq \emptyset$, which finishes the proof of theorem \ref{thm:2}.

\section{An application of theorem \ref{thm:2}}

Let $S \subset ]0, + \infty[$ be infinite and countable such that the
reverse linear ordering $>$ on $\R$ induces a well-ordering on $S$. We
saw that $\Ur$ is then indivisible but that there is no big
Ramsey degree for any $\m{X} \in \U$ as soon as $|\m{X}| \geqslant
2$. In other words, in the present context, the analogue of infinite
Ramsey's theorem holds in dimension $1$ but fails for higher
dimensions. Still, one may ask if some partition result fitting in
between holds. For example, given any $f : \funct{\Ur}{\omega}$, is
there an isometric copy of $\Ur$ inside $\Ur$ on which $f$ is
constant or injective ? It turns out that the answer is no. To see
that, consider a family $(b_n)_{n \in \omega}$ of disjoint balls
covering $\Ur$ whose sequence of corresponding radii $(r_n)_{n \in
  \omega}$ decreases strictly to $0$ and define $f :
\funct{\Ur}{\omega}$ by $f(x) = n$ iff $x \in b_n$. Then $f$ is not constant or injective on any isometric copy of $\Ur$. Observe in fact that $f$ is neither uniformly continuous nor injective on any isometric copy of
$\Ur$. However, if ``uniformly continuous'' is replaced by
``continuous'', then the result becomes true. The purpose of this
section is to provide a proof of that fact in the general case. The reader will notice the similarities with the proof of theorem \ref{thm:2}.

\begin{defn}

Let $f : \funct{\Ur}{\omega}$, $Y \subset \Ur$ and $b \in \B$ with radius $r > 0$. Say that \emph{$f$ has almost finite range on $b$ with respect to $Y$} when there is a finite family $(c_i)_{i < n}$ of elements of $\B$ with radius $r^-$ such that $f$ has finite range on $Y \cap (b \smallsetminus \bigcup _{i<n} c_i)$.  

\end{defn}

\begin{lemma}

\label{lemma:4}

Let $f : \funct{\Ur}{\omega}$ and $Y \subset \Ur$ such that for every $b \in \B$ meeting $Y$, $f$ does not have almost finite range on $b$ with respect to $Y$. Then there is an isometric copy of $\Ur$ included in $Y$ on which $f$ is injective. 

\end{lemma}

\begin{proof}

Let $\{x_n : n \in \omega \}$ be an enumeration of the elements of $\Ur$. Our goal is to construct inductively a sequence $(y_n)_{n \in \omega}$ of elements of $Y$ on which $f$ is injective and such that \[ \forall m, n \in \omega \ \ d(y_m , y_n) = d(x_m , x_n) \] 

For $y_0$, take any element in $Y$. In general, if $(y_n)_{n \leqslant k}$ is built, construct $y_{k+1}$ as follows. Consider the set $E$ defined as 
\[ E = \{ y \in \Ur : \forall \ n \leqslant k \ \ d(y,y_n) = d(x _{k+1} , x_n) \} \]

As in lemma \ref{lemma:1}, there is $b \in \B$ with radius $r > 0$ intersecting $Y$ and a set $M$ such that \[ E = b \smallsetminus \bigcup _{m \in M} B(y_m, r^-)\] 

Since $f$ does not have almost finite range on $b$ with respect to $Y$, $f$ takes infinitely many values on $E$ and we can choose $y_{k+1} \in E$ such that \[ \forall n \leqslant k \ \ f(y_n) \neq f(y_{k+1}) \ \  \] 

\end{proof}

We now turn to a proof of theorem \ref{thm:5}. Here, our strategy is
to define recursively a sequence $(Q _{\alpha})_{\alpha \in \omega
  _1}$ whose purpose is to get rid of all those parts of $\Ur$ on which $f$ is essentially of finite range: 

\begin{itemize}

\item $Q_0 = \Ur$.

\item $Q_{\alpha + 1} = Q_{\alpha} \smallsetminus \bigcup \{ b : \mathrm{f \ has \ almost \ finite \ range \ on} \ b \ \mathrm{with \ respect \ to} \ Q_{\alpha} \}$. 

\item For $\alpha < \omega _1$ limit, $Q_{\alpha} = \bigcap _{\eta < \alpha} Q_{\eta}$. 

\end{itemize}

$\Ur$ being countable, the sequence is eventually constant. Set \begin{center} $\beta = \min \{ \alpha < \omega _1 : Q_{\alpha +1} = Q_{\alpha} \}$ \end{center}

If $Q_{\beta}$ is non-empty, then $f$ and $Q_{\beta}$ satisfy the hypotheses of lemma \ref{lemma:4}. Indeed, suppose that $b \in \B$ is such that $f$ has almost finite range on $b$ with respect to $Q_{\beta}$. Then $Q_{\beta + 1} \cap b = \emptyset$. But $Q_{\beta + 1} = Q_{\beta}$ so $Q_{\beta} \cap b = \emptyset$. 

Consequently, suppose that $Q_{\beta} = \emptyset$. The intuition is that on any ball $b$, $f$ is essentially of finite range. Consequently, we should be able to show that there is $X \in \binom{\Ur}{\Ur}$ on which $f$ is continuous. 

For $\alpha < \omega _1$,
let $\mathcal{C}_{\alpha}$ be the set of all minimal elements of the collection $ \{ b : \mathrm{f \ has \ almost \ finite \ range \ on} \ b \ \mathrm{with \ respect \ to} \ Q_{\alpha} \}$. Then  

\[ \forall \ 0 < \alpha < \omega _1 \ \ Q_{\alpha} = \Ur \smallsetminus
\bigcup_{\eta < \alpha} \bigcup \mathcal{C}_{\eta} \ \ \ \ (**)\] 

\begin{claim}
Let $\alpha < \omega _1$ and $b \in \mathcal{C}_{\alpha}$. Then there is $\tilde{b} \in \binom{b}{b}$ on which $f$ is continuous.  
\end{claim}

\begin{proof}
We proceed by induction on $\alpha < \omega _1$. 

For $\alpha = 0$, let $b \in \mathcal{C}_0$. $f$ has almost finite range on $b$ with respect to $Q_0 = \Ur$ so find $c_0, \ldots , c_{n-1}$ with radius $r^-$ such that $f$ has finite range on $b \smallsetminus \bigcup_{i<n}c_i$. Then $b
\smallsetminus \bigcup_{i<n}c_i$ is isometric to $b$. Now, by theorem \ref{thm:2}, $b$ is indivisible. Therefore, there is $\tilde{b} \in \binom{b}{b}$ on which $f$ is constant, hence continuous.

Suppose now that the claim is true for every $\eta < \alpha$. Let $b \in \mathcal{C}_{\alpha}$ with radius $r \in S$. 
Find $c_0 \ldots c_{n-1} \in \B$ with radius $r^-$ and included in $b$ such that $f$ has finite range on $Q_{\alpha} \cap (b \smallsetminus \bigcup _{i<n} c_i) $. Then $b' := b \smallsetminus \bigcup _{i<n} c_i$ is isometric to $b$ and thanks to $(**)$,

\[ b' = (b'\cap Q_{\alpha}) \cup (b' \cap \bigcup _{\eta < \alpha} \bigcup \mathcal{C}_{\eta})
\]

For the same reason as in section 3, if $\mathcal{D}_{\alpha}$ is
defined as the set of all minimal elements of the collection \[ \{ c \in \bigcup _{\eta < \alpha} \mathcal{C}_{\eta} : c \subset b \wedge \forall i < n \ \ c \cap c_i = \emptyset \} \]

then we have \[ b' = (b' \cap Q_{\alpha}) \cup \bigcup \mathcal{D}_{\alpha} \]

Thanks to theorem \ref{thm:2}, $b' \cap Q_{\alpha}$ or $\bigcup \mathcal{D}_{\alpha}$ includes an isometric copy $\tilde{b}$
of $b$. If $b' \cap Q_{\alpha}$ does, then for every $i<n$, $c_i \cap
\tilde{b}$ is a metric ball of $\tilde{b}$ of same radius as
$c_i$. Thus, $\tilde{b} \smallsetminus \bigcup _{i<n} c_i$ is an
isometric copy of $b$ on which $f$ takes only finitely many values and
theorem \ref{thm:2} allows to conclude. Otherwise, suppose
that $\bigcup \mathcal{D}_{\alpha}$ includes an isometric copy of $b$. Note
that $\bigcup \mathcal{D} _{\alpha}$ includes an isometric
copy of itself on which $f$ is continuous. Indeed, by induction
hypothesis, for every $c \in \mathcal{D} _{\alpha}$, there is an
isometry $\varphi _c : \funct{c}{c}$ such that $f$ is
continuous on the range $\varphi _c ''c$ of $\varphi _c$. As in the
previous section, $\varphi := \funct{\bigcup \mathcal{D}_{\alpha}}{\bigcup \mathcal{D}_{\alpha}}$ defined as \[ \varphi = \bigcup_{c \in
  \mathcal{D}_{\alpha}} \varphi _c \] 

is an isometry. Thus, its range $\varphi '' \bigcup \mathcal{D}
_{\alpha} $ is an isometric copy of $\bigcup
\mathcal{D}_{\alpha}$ on which $f$ is continuous. Now, since
$\bigcup \mathcal{D}_{\alpha}$ includes an isometric copy of $b$, so does $\varphi '' \bigcup \mathcal{D}
_{\alpha} $ and we are done.

\end{proof}

We conclude with the same argument we used at the end of theorem \ref{thm:2}: As a metric ball, $\Ur$ is in $\mathcal{C}_{\beta}$. Thus, there is an isometric copy $X$ of $\Ur$ inside $\Ur$ on which $f$ is continuous.

\section{Concluding remarks}

The Ramsey theory of classes of finite metric spaces and their corresponding ultrahomogeneous objects is far from being fully developed so there is still a lot to investigate in that field. However, the deep connection between ultrametric spaces and trees which is exploited in this article considerably simplifies the combinatorial core of the problem and brings no help out of this specific context. Consequently, new results in the area, for example concerning Euclidean metric spaces, will presumably require the introduction of new techniques. So far the situation is not clear. We hope it will become soon.

\end{document}